\documentclass[5p,twocolumn]{elsarticle}

\usepackage[colorlinks,linkcolor=red,anchorcolor=blue,citecolor=green,CJKbookmarks=True]{hyperref}

\usepackage{graphicx}      
\usepackage{flushend}
\usepackage{natbib}        

\biboptions{authoryear}

\usepackage{amsmath,amssymb,amsfonts}
\usepackage{algorithm}
\usepackage{algorithmic}

\usepackage{savesym}
\savesymbol{AND}

\usepackage{textcomp}
\usepackage{url}
\usepackage{color}
\usepackage{tikz}
\usetikzlibrary{arrows,shapes,chains}  	
\usepackage{xcolor}
\usepackage{lipsum}
\usepackage{epstopdf}

\usepackage{amsopn}

\usepackage{enumerate}
\usepackage{cancel}
\usepackage{mathrsfs}
\usepackage{mathdots}
\usepackage{euscript}
\usepackage{amscd}
\usepackage{placeins}

\usepackage{multirow} 
\usepackage{xcolor}

\usepackage{booktabs}

\newtheorem{proposition}{Proposition}

\newtheorem{remark}{Remark}

\newdefinition{rmk}{Remark}
\newproof{pf}{Proof}

\newcommand{\bmat}{\left[ \begin{matrix}}
	\newcommand{\emat}{\end{matrix} \right]}
\newcommand{\innerprod}[2]{\langle{#1},\,{#2}\rangle}
\DeclareMathOperator{\trace}{tr}

\DeclareMathOperator{\argmax}{argmax}

\DeclareMathOperator{\E}{{\mathbb E}}
\newcommand{\Rbb}{\mathbb R}

\newcommand{\Cbb}{\mathbb C}

\newcommand{\Nbb}{\mathbb N}

\newcommand{\Zbb}{\mathbb Z}

\newcommand{\Tbb}{\mathbb T}

\newcommand{\zb}{\mathbf  z}

\newcommand{\ab}{\mathbf a}

\newcommand{\pb}{\mathbf  p}
\newcommand{\nb}{\mathbf  n}

\newcommand{\qb}{\mathbf q}
\newcommand{\tb}{\mathbf t}
\newcommand{\kb}{\mathbf k}
\newcommand{\lb}{\boldsymbol{\ell}}

\newcommand{\zerob}{\mathbf 0}

\newcommand{\Ab}{\mathbf A}
\newcommand{\Bb}{\mathbf B}

\newcommand{\Nb}{\mathbf N}

\newcommand{\Qb}{\mathbf Q}

\newcommand{\Tb}{\mathbf T}

\newcommand{\thetab}{\boldsymbol{\theta}}
\newcommand{\omegab}{\boldsymbol{\omega}}

\newcommand{\Sigmab}{\boldsymbol{\Sigma}}

\newcommand{\Lscr}{\mathscr{L}}

\newcommand{\Hcal}{\mathcal{H}}

\renewcommand{\d}{\mathrm{d}}

\newcommand{\Exp}{\mathrm{Exp}}
\newcommand{\F}{\mathrm{F}}

\newcommand{\half}{\mathrm{half}}

\begin{document}

\begin{frontmatter}

\title{A Fast Robust Numerical Continuation Solver to a Two-Dimensional Spectral Estimation Problem\tnoteref{footnoteinfo}} 
\tnotetext[footnoteinfo]{This work was supported in part by the ``Hundred-Talent Program'' of Sun Yat-sen University and the National Natural Science Foundation of China under the grant number 62103453. Conflict of interest - none declared. 
}


\author[a]{Bin Zhu\corref{cor1}}
\ead{zhub26@mail.sysu.edu.cn}
\cortext[cor1]{Corresponding author}
\author[a]{Jiahao Liu}
\ead{liujh226@mail2.sysu.edu.cn}
\address[a]{School of Intelligent Systems Engineering, Sun Yat-sen University, Waihuan East Road 132, 510006 Guangzhou, China}




\begin{abstract}
This paper presents a fast algorithm to solve a spectral estimation problem for two-dimensional random fields. The latter is formulated as a convex optimization problem with the Itakura-Saito pseudodistance as the objective function subject to the constraints of moment equations. We exploit the structure of the Hessian of the dual objective function in order to make possible a fast Newton solver. Then we incorporate the Newton solver to a predictor-corrector numerical continuation method which is able to produce a parametrized family of solutions to the moment equations. We have performed two sets of numerical simulations to test our algorithm and spectral estimator. The simulations on the frequency estimation problem shows that our spectral estimator outperforms the classical windowed periodograms in the case of two hidden frequencies and has a higher resolution. The other set of simulations on system identification indicates that the numerical continuation method is more robust than Newton's method alone in ill-conditioned instances.
\end{abstract}


\end{frontmatter}

\section{Introduction}
Spectral estimation is a classical problem in signal processing 
closely related to the field of systems and control, as it finds applications in e.g., stochastic realization, system identification and modeling, estimation and filtering \citep{stoica2005spectral,LP15}.
The problem is about estimating the power spectral density function, which describes the statistical power distribution over the frequency domain, of a zero-mean second-order stationary random field from a finite number of measurements.
Once such a spectral density is reconstructed, one can then perform \emph{spectral factorization} to obtain a shaping filter, i.e., a dynamical system. When the filter is excited with white noise, it is able to reproduce a random field that is statistically close to the original one that gives the measurements.

A still active line of research on the spectral problem follows the idea of \emph{rational covariance extension} which was initially formulated in \citet{Kalman} and subsequently developed in \citet{Gthesis,Georgiou-87,BLGM-95,byrnes1997partial}.
The central object of investigation is a set of linear integral equations, called \emph{trigonometric moment equations} \citep{Grenander_Szego,A65moments,KreinNudelman}, for which we want to find \emph{rational} solutions.
The additional requirement of rationality is natural from the systems-theoretic viewpoint, as it ensures that the corresponding shaping filter is linear-time invariant (LTI) and can be realized as a finite-dimensional system.
Typically, the number of integral equations is finite and there are usually infinitely many solutions which are not necessarily rational.
In order to cure such ill-posedness of the problem and to promote rational solutions, convex optimization was then incorporated into the framework. 
More precisely, moment equations are treated as equality constraints, and one searches for a candidate solution that maximizes a certain entropy functional  \citep{BGL-98,BGL-01,byrnes2001cepstral} or minimizes a suitable pseudodistance from a \emph{prior} spectral density \citep{Georgiou-L-03,FPR-07,FPR-08,RFP-09,RFP-10-wellposedness,FMP-12,Z-14rat,Z-14,Z-15}. 
We mention in particular, the seminal work in \cite{BGL-THREE-00} where a flexible filter-bank was introduced and the resulting spectral estimator, called ``THREE'' (Tunable High REsolution Estimator), exhibits high-resolution properties.

The aforementioned theories are established for random processes and one-dimensional systems, i.e., random fields and dynamical systems that depend on one index, in most cases the time.	
Motivated by many practical applications involving multidimensional systems and random fields such as image processing \citep{ekstrom1984digital} and parameter estimation in automotive radar systems \citep{rohling2012continuous,engels2014target,ZFKZ2019fusion}, the research in rational covariance extension has also been extended to the multidimensional case, see \citet{Georgiou-06,KLR-16multidimensional,RKL-16multidimensional,ringh2018multidimensional,ZFKZ2019M2,ZFKZ2020M2-SIAM,Zhu-M2-LineSpec,Zhu-Zorzi-2021-cepstral}. 
Discrete versions of the theory that facilitate numerical computation have also been developed in \citet{ringh2015multidimensional,Zhu-Zorzi-SYSID21}, following the idea in \citet{CFPP-11,LPcirculant-13} for the 1-d case.

However, we want to point out that while theoretical developments in this area are significant, algorithmic studies seem scarce. We mention the works \citet{FRT-11,baggio2018further} on an integral-form iterative algorithm, \citet{carli2013efficient} on block-Toeplitz matrix completion, \citet{ringh2015afast} on a fast Newton solver, and \citet{enqvist2001homotopy,zhu2018wellposed} on numerical continuation methods, all in the 1-d case. As a first step towards efficient algorithms in the multidimensional case, in this paper we propose a fast implementation of the classical Newton's method to solve a 2-d spectral estimation problem formulated as a moment-constrained optimization problem where the \emph{Itakura-Saito} pseudodistance is used as the objective function.
Such a choice of the objective function first appeared in \cite{enqvist2008minimal} and later further developed in \citet{FMP-12,ZFKZ2019M2} where it was shown that the solution is indeed rational and of bounded complexity (in terms of the McMillan degree).
For the convenience of numerical computation, we mainly deal with the dual optimization problem which typically has much fewer number of variables than the primal problem. 
In order to achieve efficiency in inverting the Hessian of the dual objective function, which is well known to be the major computational burden of Newton's method, we exploit the special Toeplitz-block Toeplitz structure of the Hessian matrix. 	
Then we integrate the fast Newton solver to a more sophisticated numerical continuation method \citep{allgower2003introduction} which brings more robustness when the solution is close to the boundary of the feasible set (see Subsec.~\ref{subsec: ill_cond_ex}). 
Finally, we apply our algorithm to solve a 2-d frequency estimation problem. Simulation results show that our spectral estimator performs favorably compared with traditional periodogram-based techniques in terms of both cumulative estimation errors and frequency resolution.

The outline of the paper is as follows: Section \ref{Sec:problem statement} provides an optimization approach for solving a 2-d spectral estimation problem. In Section \ref{Sec:The structure of the Hessian} we show that the Hessian has a Toeplitz-block Toeplitz structure which can be used to implement a fast inversion algorithm. 
In Section \ref{Sec:A numerical continuation solver} we develop a numerical continuation solver for our optimization problem and show its robustness in an ill-conditioned example.
Next in Sections \ref{Sec:Fre_est} and \ref{sec:Sysid}, we apply our approach to problems of frequency estimation and model approximation in system identification,
respectively. Finally, Section \ref{Sec:Conclusions} draws the conclusions.

\section{Problem formulation}\label{Sec:problem statement}
Consider a (scalar) complex-valued zero-mean second-order stationary $2$-d random field $\{y(\tb),\tb \in \Zbb^2\}$. 
By stationarity, we mean that the covariance lags  $\sigma_{\kb}:=\E y(\tb+\kb)\overline{y(\tb)}$ depend on the vector of difference $\mathbf{k}$ only.
In many scientific and engineering fields, we face the problem of estimating the power spectral density of some underlying stationary random field from a finite number of measurements 
\begin{equation}
	\{ y(\mathbf{t}) : 0\leq t_j\leq T_j-1,\ j=1,2\}.
\end{equation}
Here the power spectral density of the random field is defined as the discrete-time Fourier transform (DTFT) of the covariance lags:
\begin{equation}\label{1}
	\Phi(e^{i\boldsymbol{\theta}}) := \sum_{\mathbf{k}\in \Zbb^2}\sigma_\mathbf{k}e^{-i\langle\mathbf{k},\boldsymbol{\theta}\rangle},
\end{equation}
where $\thetab=[\theta_1,\theta_2]\in \Tbb^2=\left[0,2\pi\right)^2$, $e^{i\boldsymbol{\theta}}$ is a shorthand for $(e^{i\theta_1},e^{i\theta_2})$, and $\langle \mathbf{k},\boldsymbol{\theta} \rangle=k_1\theta_1+k_2\theta_2$ is the inner product in $\mathbb{R}^2$.
The spectral density $\Phi$ is viewed as a function of the frequency vector $\thetab$ but we keep the conventional notation $\Phi(e^{i\thetab})$.
Conversely, the covariances $\sigma_\mathbf{k} $ are the Fourier coefficients of the power spectral density $\Phi(e^{i\thetab})$, that is,
\begin{equation}\label{2}
	\sigma_\mathbf{k}=\int_{\mathbb{T}^2}e^{i\langle \mathbf{k},\boldsymbol{\theta} \rangle}\Phi(e^{i\boldsymbol{\theta}})\mathrm{d}m(\boldsymbol{\theta}),
\end{equation}
where $\mathrm{d}m\left(\boldsymbol{\theta}\right)=\frac{1}{{(2\pi)}^2}\mathrm{d}\theta_1\d \theta_2$ is the normalized Lebesgue measure in $\mathbb{T}^2$. 
Later on, the variable $\thetab$ or $e^{i\thetab}$ of integration is usually omitted when it is clear from the context.
Equation \eqref{2} is also called a \emph{trigonometric moment equation}, and it forms the basis of the spectral estimation approach via covariance extension.
The main steps of the latter approach are outlined as follows:
\begin{enumerate}
	\item Fix an index set $\Lambda$, e.g.,
	\begin{equation}\label{set_Lambda}
		\Lambda = \{\kb=(k_1,k_2)\in\Zbb^2 : |k_j| \leq n_j, j=1,2\},
	\end{equation}
	where each positive integer $n_j$ satisfies $n_j\ll T_j$.
	\item Estimate the covariance lags  $\sigma_{\kb}$ for $\kb \in \Lambda$ by the standard time average:
	\begin{equation}\label{cov_est_ave}
		\hat{\sigma}_{\kb} = \frac{1}{T_1T_2} \sum_{\tb} y(\tb+\kb) \overline{y(\tb)},
	\end{equation}
	where the summation index $\tb$ is chosen such that the summands are well defined.
	\item Solve the integral equations \eqref{2} indexed by $\kb\in\Lambda$ with $\sigma_{\kb}$ replaced by the estimate $\hat{\sigma}_{\kb}$ for a spectral density $\Phi$.
\end{enumerate}
In other words, we aim to find a spectral density function $\Phi$ that matches the estimated covariance lags $\{\hat{\sigma}_{\kb}\,:\,\kb\in\Lambda\}$. In the following, we omit the hat in the covariance data for the notational convenience and focus on the solution to the moment equations.

Typically, the integral equations admit infinitely many solutions. In order to promote uniqueness of the solution,  a common practice in the literature is built on convex optimization. In particular, we take the $2$-d version of the Itakura-Saito (IS) pseudo-distance as the objective function and treat the moment equations as equality constraints. In this way the spectral estimation problem is formulated as the following optimization problem:
\begin{equation}\label{optimization problem}
	\begin{aligned}
		&
		\begin{aligned}
			\min_{\Phi > 0} \ D(\Phi,\Psi) := \int_{\mathbb{T}^2} \left[ \log (\Phi^{-1}\Psi) + \Psi^{-1}(\Phi-\Psi) \right] \mathrm{d}m
		\end{aligned}
		\\
		& \text{s.t.} \ \sigma_\mathbf{k}=\int_{\mathbb{T}^2}e^{i\langle \mathbf{k},\boldsymbol{\theta} \rangle}\Phi\, \mathrm{d}m
		 \quad  \forall \kb \in \Lambda,
	\end{aligned}
\end{equation}
where the spectral density $\Phi$ is the optimization ``variable'' and $\Psi$ is called a \emph{prior} which represents an extra piece of data that embeds some \emph{a priori} information on the desired solution. In the absence of such information, $\Psi$ can be chosen as a constant which corresponds to a white noise.
The IS pseudo-distance $D\left(\Phi,\Psi\right)$ has the properties that  $D\left(\Phi,\Psi\right)\geq0$, and the equality holds if and only if $\Phi=\Psi$ almost everywhere, see e.g., \cite{LP15}. Moreover, the solution to the above optimization problem turns out to be \emph{rational}, and enjoys a low-complexity property  in the 1-d case  as discussed in \cite{FMP-12}. Rationality is important from the systems-theoretic viewpoint because it is directly connected to the shaping filter that can regenerate the original random field (in a statistical sense when fed with white noise) and the physical realizability of such a filter.

In what follows, we mainly work on the discrete version of the optimization problem in \eqref{optimization problem} for two reasons.
First, the implementation of any numerical algorithm on a computer to solve the optimization problem ultimately involves discretization. Hence, we may also discretize the problem in the first place.
Second, a discrete spectrum admits an interpretation as a periodic random field \citep{ZFKZ2019M2}, which is a practical model for the collected \emph{finite} measurements, provided that the period is sufficiently large in each dimension.

We need to introduce a number of symbols in order to state the discrete optimization problem.
Define first a vector $\Nb:=[N_1, N_2]\in\Nbb^2_+$ and then a finite index set 
\begin{equation}
	\mathbb{Z}_{\mathbf{N}}^2:=\left\{
	\lb=(\ell_1,\ell_2):0 \leq \ell_j \leq N_j-1,j=1,2\right\},
\end{equation}
whose cardinality is equal to $|\Nb|:=N_1N_2$.
Then the discretization of $\mathbb{T}^2$ can be written as
\begin{equation}
	\mathbb{T}_\mathbf{N}^2:=\left\{\left( \frac{2\pi}{N_1}\ell_1,\frac{2\pi}{N_2}\ell_2\right):\lb \in \mathbb{Z}_\mathbf{N}^2 \right\}.
\end{equation}
Moreover, let $\boldsymbol{\zeta}_{\lb}:=(\zeta_{\ell_1},\zeta_{\ell_2})$ be the element of the discretized 2-torus with $\zeta_{\ell_j}=e^{i2\pi \ell_j/N_j}$.
A discrete measure with equal mass $1/|\Nb|$ on the grid points in $\Tbb^2_\Nb$ is defined as
\begin{equation}
	\d \nu_{\mathbf{N}}(\thetab) := \sum_{\boldsymbol{\lb} \in \mathbb{Z}_{\mathbf{N}}^2}\delta (\theta_1-\frac{2\pi}{N_1}\ell_1,\theta_2-\frac{2\pi}{N_2}\ell_2) \frac{\d \theta_1 \d \theta_2}{N_1 N_2},
\end{equation}
where $\delta(\cdot)$ is the Dirac delta.
Thus integrals against $\d\nu_{\Nb}$ is understood as (normalized) Riemann sums:
\begin{equation}
	\int_{\mathbb{T}^2}f(e^{i\thetab})\mathrm{d} \nu_\mathbf{N} = \frac{1}{|\mathbf{N}|}\sum_{\boldsymbol{\lb} \in \mathbb{Z}_\mathbf{N}^2}f(\boldsymbol{\zeta}_{\boldsymbol{\lb}}).
\end{equation}
Then the discrete version of the IS distance takes the form
\begin{equation}
	D_{\mathbf{N}}(\Phi,\Psi) := \int_{\mathbb{T}^2} \left[ \log (\Phi^{-1}\Psi) + \Psi^{-1}(\Phi-\Psi) \right] \mathrm{d}\nu_\mathbf{N},
\end{equation}
and the discrete optimization problem is
\begin{equation}\label{discretized optimization problem}
	\begin{aligned}
		&\min _{\Phi(\boldsymbol{\zeta}_{\lb})> 0,\ \forall \lb \in \mathbb{Z}_\mathbf{N}^2} \ D_{\mathbf{N}}(\Phi, \Psi)\\
		&\text{s.t.} \  \sigma_\mathbf{k}=\int_{\mathbb{T}^2}e^{i\langle \mathbf{k},\boldsymbol{\theta} \rangle} \Phi\, \d \nu_\mathbf{N} \quad  \forall \kb \in \Lambda.
	\end{aligned}
\end{equation}

According to \cite{ZFKZ2019M2}, the optimization problem above is well-posed, and the optimal spectral density, defined on the discrete grid $\Tbb^2_{\Nb}$, has the form
\begin{equation}\label{optimal solution}
	\hat{\Phi}=(\Psi^{-1} + \hat{Q})^{-1},
\end{equation}
where $\hat{Q}(e^{i\thetab}):=\sum_{\kb \in \Lambda} \hat{q}_{\kb} e^{-i\innerprod{\kb}{\thetab}}$ is a trigonometric Laurent polynomial corresponding to the optimal solution of the (equivalent) dual problem:
\begin{equation}\label{J_dual}
	\underset{\Qb\in\Lscr_+}{\text{min}} \ J_\Psi(\Qb):=\langle\Qb,\Sigmab\rangle-\int_{\Tbb^2}\log(\Psi^{-1}+Q)\d\nu_\Nb.
\end{equation}
The notation is briefly explained next.
\begin{itemize}
	\item The variable $\Qb=\{q_\kb\}_{\kb\in\Lambda}$ contains the Lagrange multipliers such that each $q_\kb\in\Cbb$, $q_{-\kb} = \overline{q_\kb}$, and $q_\zerob$ is real where the index set $\Lambda$ can be taken as \eqref{set_Lambda}. 
	\item $\Sigmab=\{\sigma_\kb\}_{\kb\in\Lambda}$ consists of the covariance data of the underlying random field.
	\item $\innerprod{\Qb}{\Sigmab}:=\sum_{\kb\in\Lambda} q_\kb \overline{\sigma_{\kb}}$ is an inner product between multisequences indexed in $\Lambda$ which is real-valued due to the symmetry in $\Qb$ and $\Sigmab$ with respect to the origin.
	\item The feasible set
	\begin{equation}\label{feasible_set}
		\Lscr_+:=\left\{ \Qb : (\Psi^{-1}+Q)(\thetab)>0,\ \forall \thetab\in\Tbb^2_\Nb \right\}
	\end{equation}
	so that the logarithm inside the integral is well defined.
\end{itemize}

We shall take the grid size $\Nb$ sufficiently large so that Assumption 2 in \citet{ZFKZ2019M2} is satisfied, leading to the strict convexity of the dual problem \eqref{J_dual}.

\section{Structure of the Hessian}\label{Sec:The structure of the Hessian}
In what follows, we are mainly concerned about the numerical solution of the optimization problem \eqref{J_dual} via Newton's method, for which we need to compute the gradient and the Hessian of the dual objective function $J_\Psi(\Qb)$. Moreover, we hope to exploit the structure of the Hessian in order to boost the solution speed, notably in the computation of the Newton direction. To this end, we need first to fix an orthogonal basis for the vector space in which the dual variable $\Qb$ resides. In the specific formulation above, the dual variable can be identified as a $(2n_2+1)\times(2n_1+1)$ complex matrix
\begin{equation}\label{mat_Q}
	\Qb = \bmat
	q_{-n_1,-n_2} & \cdots & q_{0,-n_2} & \cdots & q_{n_1,-n_2} \\
	\vdots & & \vdots & & \vdots \\
	q_{-n_1,0} & \cdots & q_{0,0} & \cdots & q_{n_1,0} \\
	\vdots & & \vdots & & \vdots \\
	q_{-n_1,n_2} & \cdots & q_{0,n_2} & \cdots & q_{n_1,n_2} \\
	\emat.
\end{equation}
Using this notation, the inner product in the dual function can be computed as 
\begin{equation}\label{inner_prod_mat}
	\innerprod{\Qb}{\Sigmab} = \trace(\Qb\Sigmab^*)
\end{equation}
where the matrix $\Sigmab$ is defined similarly. Moreover, let us also introduce the matrix of $2$-d complex exponentials
\begin{equation}
	\begin{aligned}
		\Exp(\thetab) & = \bmat
		e^{i(-n_1\theta_1 -n_2\theta_2)} & \cdots & e^{i(n_1\theta_1 -n_2\theta_2)} \\
		\vdots & & \vdots \\
		e^{i(-n_1\theta_1 +n_2\theta_2)} & \cdots & e^{i(n_1\theta_1 +n_2\theta_2)} \\
		\emat \\
		& = \ab_{n_2}(\theta_2) \ab_{n_1}(\theta_1)^\top
	\end{aligned}
\end{equation}
where
\begin{equation}
	\ab_n(\theta) = \bmat e^{-in\theta} & \cdots & 0 & \cdots & e^{in\theta}\emat^\top
\end{equation}
is a column vector of dimension $2n+1$. Then the polynomial evaluation can be written as
\begin{equation}
	Q(e^{i\thetab}) = \trace(\Qb\,\Exp(\thetab)^*) = \innerprod{\Qb}{\Exp(\thetab)}.
\end{equation}

For each $\kb \in\Lambda$, define a matrix $E_{\kb}$ which is composed of $\{q_{\lb}\}_{\lb\in\Lambda}$ in the way of \eqref{mat_Q} such that $q_\kb=1$ while the rest $q_{\lb}$ all equal to zero. As a result, we have for any $\Qb$ the expansion
\begin{equation}\label{expansion_Q}
	\Qb = \sum_{\kb\in\Lambda} q_\kb E_\kb.
\end{equation}
Obviously, the set of matrices $\{E_\kb\}_{\kb\in\Lambda}$ are orthogonal with respect to the inner product \eqref{inner_prod_mat}, and they form a basis for the ambient vector space of $\Qb$ over the field of \emph{complex} numbers.

\begin{remark}
	The matrix $\Qb$ in \eqref{mat_Q} lives in an ambient space of complex dimension $M:=(2n_1+1)(2n_2+1)$ which is isomorphic to $\Rbb^{2M}$. However, with the additional symmetry $q_{-\kb} = \overline{q_{\kb}}$, the matrix $\Qb$ actually belongs to a subspace of real dimension $M$. To count the dimension, just notice the independent variables
	\begin{itemize}
		\item $q_{0,0}$, $1$ real variable,
		\item $q_{0,1},\dots,q_{0,n_2}$, $n_2$ complex variables,
		\item \begin{equation*}
			\begin{matrix}
				q_{1,-n_2} & \cdots & q_{1,0} & \cdots & q_{1,n_2} \\
				\vdots & & \vdots & & \vdots \\
				q_{n_1,-n_2} & \cdots & q_{n_1,0} & \dots & q_{n_1,n_2} \\
			\end{matrix},
		\end{equation*}
		$n_1(2n_2+1)$ complex variables.
	\end{itemize}
	In that case, it is also observed that $\Qb$ is invariant if it undergoes a rotation by $180\deg$ followed by complex conjugation.
	Although it may be tempting to write down an expansion similarly to \eqref{expansion_Q} directly in a basis of the subspace with \emph{real} coefficients, as we will show next, working with complex coordinates in fact results in a simpler structure of the Hessian of $J_\Psi(\Qb)$.
\end{remark}

Before doing the computation in the basis $\{E_\kb\}_{\kb\in\Lambda}$, let us define a linear operator
\begin{equation}
	\Gamma:\, \Phi \mapsto \left\{ \int_{\Tbb^2} e^{i\innerprod{\kb}{\thetab}} \Phi\, \d\nu_\Nb \right\}_{\kb\in\Lambda}
\end{equation}
that sends a (discrete) spectral density to its Fourier coefficients indexed in $\Lambda$ which can then be organized into a matrix consistent with \eqref{mat_Q}. According to \cite{ZFKZ2019M2}, the first-order differential of $J_\Psi$ can be written compactly as
\begin{align}
	\delta J_\Psi (\Qb;\delta\Qb) = \innerprod{\delta\Qb}{\Sigmab -\Gamma\left( (\Psi^{-1} +Q)^{-1} \right)}.
\end{align}
It then follows from standard results in Calculus that
\begin{equation}\label{partial_J_dual}
	\begin{aligned}
		\frac{\partial J_\Psi(\Qb)}{\partial \overline{q_\kb}} & = \delta J_\Psi (\Qb;E_{-\kb}) \\
		& = \sigma_{\kb} - \int_{\Tbb^2} e^{i\innerprod{\kb}{\thetab}} (\Psi^{-1} +Q)^{-1} \d\nu_\Nb,
	\end{aligned}
\end{equation}
which is called the \emph{Wirtinger derivative} \citep[see e.g.,][]{kreutz2009complex} where $q_\kb$ and $\overline{q_{\kb}}$ are treated as ``independent'' variables in the formal calculation. The \emph{gradient} of $J_\Psi$ is defined as $\nabla J_\Psi(\Qb):=$
\begin{equation}\label{gradient_J}
	\begin{aligned}
		\left[
		\begin{array}{ccccc}
			\dfrac{\partial J_\Psi(\Qb)}{\partial \overline{q_{0,0}}} 
			& \dfrac{\partial J_\Psi(\Qb)}{\partial \overline{q_{0,1}}} & \cdots 
			& \dfrac{\partial J_\Psi(\Qb)}{\partial \overline{q_{0,n_2}}} & \dfrac{\partial J_\Psi(\Qb)}{\partial \overline{q_{1,-n_2}}}
		\end{array}
		\right. \\
		\left.
		\begin{array}{cccccc}
			\cdots & \dfrac{\partial J_\Psi(\Qb)}{\partial \overline{q_{1,n_2}}} & \cdots
			& \dfrac{\partial J_\Psi(\Qb)}{\partial \overline{q_{n_1,-n_2}}} & \cdots
			& \dfrac{\partial J_\Psi(\Qb)}{\partial \overline{q_{n_1,n_2}}}
		\end{array}
		\right]^\top,
	\end{aligned}
\end{equation}
where the complex partial derivatives are organized into a column vector of size $n_2+1+n_1(2n_2+1)$ using the lexicographical ordering.

The second-order differential of $J_\Psi$ at $\Qb$ can be written as a quadratic form
\begin{equation}
	\begin{aligned}
		& \delta^2 J_\Psi(\Qb;\delta\Qb^{(1)},\delta\Qb^{(2)}) \\ = & \sum_{\kb\in\Lambda} \delta q_\kb^{(1)} \sum_{\lb\in\Lambda} \delta q^{(2)}_{\lb} \delta^2 J_\Psi(\Qb;E_\kb,E_{\lb}),
	\end{aligned}
\end{equation}
where $\{\delta q_{\kb}^{(j)}\}_{\kb\in\Lambda}$ are the coordinates of $\delta\Qb^{(j)}$ in the sense of \eqref{expansion_Q} for $j=1, 2$. Therefore, we can identify the second-order partial derivatives as
\begin{equation}
	\begin{aligned}
		\frac{\partial^2 J_\Psi(\Qb)}{\partial q_{\lb} \partial \overline{q_\kb}} & := \frac{\partial}{\partial q_{\lb}} \left[ \frac{\partial J_\Psi(\Qb)}{\partial \overline{q_\kb}} \right] \\
		& = \int_{\Tbb^2} e^{i\innerprod{\kb -\lb}{\thetab}} (\Psi^{-1}+Q)^{-2} \d\nu_\Nb
	\end{aligned}
\end{equation}
where the expression follows from \cite{ZFKZ2019M2}.
Similar to the vector notation in \eqref{gradient_J}, we can collect the second-order partials into a matrix
\begin{equation}\label{Hessian_J}
	\Hcal(\Qb) = \left[ \frac{\partial^2 J_\Psi(\Qb)}{\partial q_{\lb} \partial \overline{q_\kb}} \right]
\end{equation}
where the ``row index'' $\kb$ and the ``column index'' $\lb$ range in the set
\begin{equation}\label{Lambda_half}
	\Lambda_{\half} := \{0\}\times\{0,\dots,n_2\} \cup \{1,\dots,n_1\} \times \{-n_2,\dots,n_2\},
\end{equation}
roughly half of the index set $\Lambda$ in \eqref{set_Lambda}.
We shall refer to \eqref{Hessian_J} as the \emph{Hessian} of $J_\Psi$. It follows from the strict convexity of $J_\Psi$ that $\Hcal(\Qb)$ is (Hermitian) positive definite.
The structure of the Hessian can be made more explicit by defining the scalar quantity
\begin{equation}
	h_\kb(\Qb) = \int_{\Tbb^2} e^{i\innerprod{\kb}{\thetab}} (\Psi^{-1}+Q)^{-2} \d\nu_\Nb.
\end{equation}
For the moment, let us drop the dependence on $\Qb$. Then the Hessian can be partitioned as
\begin{equation}
	\Hcal = \bmat A & B^* \\ B & C \emat
\end{equation}
where
\begin{equation}
	A = \bmat h_{0,0} & h_{0,-1} & \cdots & h_{0,-n_2} \\
	h_{0,1} & h_{0,0} & \cdots & h_{0,-n_2+1} \\
	\vdots & \ddots & \ddots & \vdots \\
	h_{0,n_2} & \cdots & h_{0,1} & h_{0,0} \emat
\end{equation}
is $(n_2+1)\times(n_2+1)$ Toeplitz,
\begin{equation}
	B = \bmat B_1 \\ \vdots \\ B_{n_1} \emat 
\end{equation}
is an $n_1$-block vector where each block
\begin{equation}
	B_j = \bmat h_{j,-n_2} & h_{j,-n_2-1} & \cdots & h_{j,-2n_2} \\
	h_{j,-n_2+1} & h_{j,-n_2} & \cdots & h_{j,-2n_2+1} \\
	\vdots & \vdots & \ddots  & \vdots \\
	h_{j,n_2} & \cdots & h_{j,n_2-1} & h_{j,0} \emat
\end{equation}
is $(2n_2+1)\times (n_2+1)$ Toeplitz, and
\begin{equation}\label{C-matrix}
	C = \bmat H_0 & H_{-1} & \cdots & H_{-n_1+1} \\
	H_1 & H_0 & \cdots & H_{-n_1+2} \\
	\vdots & \ddots & \ddots & \vdots \\
	H_{n_1-1} & \cdots & H_{1} & H_{0} \emat
\end{equation}
is a Toeplitz-block Toeplitz (abbreviated as TBT) matrix\footnote{A TBT matrix is also known as a $2$-level Toeplitz matrix.} such that each block
\begin{equation}
	H_j = \bmat h_{j,0} & h_{j,-1} & \cdots & h_{j,-2n_2} \\
	h_{j,1} & h_{j,0} & \cdots & h_{j,-2n_2+1} \\
	\vdots & \ddots & \ddots & \vdots \\
	h_{j,2n_2} & \cdots & h_{j,1} & h_{j,0} \emat
\end{equation}
is $(2n_2+1)\times(2n_2+1)$ Toeplitz.

Newton's method for the minimization problem \eqref{J_dual} involves solving the linear system of equations
\begin{equation}\label{compute_Newton_direction}
	\Hcal(\Qb) \delta \qb = -\nabla J_\Psi(\Qb)
\end{equation}
for $\delta\qb$, and this can be done via standard block elimination. More specifically, consider the linear equations
\begin{equation}\label{block_linear_eqns}
	\bmat A & B^* \\ B & C \emat \bmat x \\ y \emat = \bmat a \\ b \emat.
\end{equation}
Assuming that $C$ is invertible, we have
\begin{equation}
	( A - B^*C^{-1}B )\, x = a - B^* C^{-1} b
\end{equation}
where the matrix on the left-hand side is known as the \emph{Schur complement} of the block $C$ of the coefficient matrix in \eqref{block_linear_eqns}. In our specific problem \eqref{compute_Newton_direction} where the Hessian is positive definite, the Schur complement is guaranteed to be positive definite as well. Consequently, the vector $x$ can be solved by inverting an $(n_2+1)\times(n_2+1)$ matrix. Then the remaining $y$ can be recovered from the second block equation in \eqref{block_linear_eqns} as
\begin{equation}
	y = C^{-1} (b-Bx).
\end{equation}
Apparently, the major computational burden is caused by the inversion of the TBT matrix $C$ of dimension $n_1(2n_2+1)$, for which an efficient algorithm can be found in \citet{wax1983efficient}.

\begin{remark}\label{footnote_quasi_Newton}
	The procedure above is in fact, a \emph{quasi-Newton} method instead of the true Newton's method which involves the \emph{full} Hessian of $J_\Psi$ with respect to all the complex variables and their conjugates in the matrix $\Qb$ in \eqref{mat_Q}. See the discussion around Eq.~(112) in \citet{kreutz2009complex}. In principle, the convergence of a quasi-Newton method could be slower than the quadratic convergence of Newton's method. 
\end{remark}

\begin{remark}
	The paper \citet{ringh2015afast} studies a fast algorithm to solve the circulant rational covariance extension problem in the scalar $1$-dimensional case. The authors of that paper considered real processes, so that the symmetry between the variables reduces to $q_{-k} = q_k$. Their Hessian has a Toeplitz-plus-Hankel structure which was also noticed in an earlier work \citep{SIGEST-01}. We expect our Hessian to have a similar structure (more than TBT) when all the second-order partial derivatives are taken into consideration. In that case, it also seems possible to devise a fast inversion algorithm for the full Hessian similar to the one in \citet{wax1983efficient}. 
\end{remark}

\section{A numerical continuation solver}\label{Sec:A numerical continuation solver}

In solving the optimization problem \eqref{J_dual} using Newton's method, it is very often noticed that as the iterates go near the boundary of the feasible set $\Lscr_+$, the condition number of the Hessian grows significantly, making it hard to achieve convergence. To handle such a situation, a more numerically stable method to compute the minimizer is available in the literature under the name ``numerical continuation'' \citep{allgower2003introduction,zhu2018wellposed}. The basic idea goes as follows.

Due to the fact that the problem \eqref{J_dual} has a unique interior solution in the domain $\Lscr_+$ \citep[cf.~][]{ZFKZ2019M2}, it is  equivalent to directly solving the stationary-point equation
\begin{equation}\label{stationary_point_eqn}
	J_2(\Psi,\Qb) = \zerob
\end{equation}
subject to the constraint $\Qb\in\Lscr_+$ given the prior $\Psi$. Note that we have changed the notation to include the explicit dependence of $J$ on $\Psi$, and $J_2\equiv\nabla_{\Qb}\, J$ denotes taking gradient with respect to the ``second'' variable $\Qb$ as indicated in \eqref{gradient_J}. Instead of solving one single equation \eqref{stationary_point_eqn}, a numerical continuation method deals with a parametrized family of such equations. More specifically, suppose that we want to solve \eqref{stationary_point_eqn} or the optimization problem \eqref{J_dual} for $\Qb_1$ given the prior $\Psi=\Psi_1$, which is hard. On the contrary, the solution $\Qb_0$ is ``easy'' to compute when one takes $\Psi=\Psi_0$. Then one hopes to depart from $\Qb_0$ and navigate to the desired $\Qb_1$ along a continuous curve. Indeed, such a curve-tracing algorithm is possible. In order to achieve it, one first constructs a continuous deformation from $\Psi_0$ to $\Psi_1$, called a \emph{homotopy}, the simplest one being the convex combination
\begin{equation}
	h(t,\thetab) := (1-t) \Psi_0(\thetab) + t\Psi_1(\thetab),
\end{equation}
where the real parameter $t\in U:=[0,1]$ and $\thetab\in\Tbb^2$. Apparently, we have $h(0,\thetab) = \Psi_0(\thetab)$ and $h(1,\thetab) = \Psi_1(\thetab)$. Moreover, $h$ is jointly continuous in $(t,\thetab)$ as long as both $\Psi_0$ and $\Psi_1$ are continuous in $\thetab$. 
For each fixed $t\in U$, $\Psi_t(\thetab) := h(t,\thetab)$ is a spectral density due to positivity. We can now write down the equation
\begin{equation}\label{parametrized_eqns}
	J_2 (h(t,\cdot),\Qb) = \zerob
\end{equation}
parametrized by $t\in U$. For each fixed $t$, we are in fact solving \eqref{stationary_point_eqn} with $\Psi=\Psi_t$, and we know that there exists a unique solution which can be written as $\Qb(t)$. As a consequence of \citet[Theorem~3]{ZFKZ2019M2}, the solution curve $\Qb(t)$ is differentiable with respect to $t$.

Next, we shall derive the ordinary differential equation (ODE) that $\Qb(t)$ satisfies.
Substitute $\Qb$ in \eqref{parametrized_eqns} with $\Qb(t)$, and differentiate both sides of the equation with respect to $t$, yielding
\begin{equation}\label{derivative_eqn}
	J_{21} (h(t,\cdot),\Qb(t); h_1(t,\cdot)) + J_{22} ((h(t,\cdot),\Qb(t); \Qb'(t))) = \zerob.
\end{equation}
The semicolon notation here means e.g.,
\begin{equation*}
	J_{21} (\,\cdot\,,\,\cdot\,; \xi) = J_{21} (\,\cdot\,,\,\cdot\,) (\xi),
\end{equation*}
the linear operator $J_{21} (\,\cdot\,,\,\cdot\,)$ applied to the object $\xi$. Clearly, such an operation reduces to a suitable matrix-vector product if $\xi$ resides in a finite-dimensional space. Now, we need to obtain an explicit expression for the operators in \eqref{derivative_eqn}. Stack the complex variables $\{q_\kb\}_{\kb\in\Lambda_{\half}}$ indexed in the ``half'' set \eqref{Lambda_half} into a long column vector $\qb$ in accordance with the variable arrangement in the Hessian \eqref{Hessian_J}. Then we have the term
\begin{equation}
	J_{22} ((h(t,\cdot),\Qb(t); \Qb'(t))) \approx \Hcal(\Psi_t,\Qb(t))\, \qb'(t),
\end{equation}
where the approximation comes from the fact that we have ignored the complex derivatives with respect to the conjugate variables indexed in $\Lambda\backslash\Lambda_{\half}$, similar to what we have done in defining the Hessian \eqref{Hessian_J}, see Remark~\ref{footnote_quasi_Newton}.
Computation of the other term in \eqref{derivative_eqn} is carried out as follows. The partial derivatives of $J_\Psi(\Qb)$ are given in \eqref{partial_J_dual}, from which we can compute the $\kb$-th component of the linear operator $J_{21}(\Psi,\Qb)$:
\begin{equation}
	\begin{aligned}
		& J_{21,\kb}(\Psi,\Qb;\delta\Psi) \\
		= & D_{\Psi} \left(\sigma_{\kb} - \int_{\Tbb^2} e^{i\innerprod{\kb}{\thetab}} (\Psi^{-1} +Q)^{-1} \d\nu_\Nb\right) \\
		= & - \int_{\Tbb^2} e^{i\innerprod{\kb}{\thetab}} (\Psi^{-1} +Q)^{-2} \Psi^{-2} \delta\Psi\, \d\nu_\Nb
	\end{aligned}
\end{equation} 
which is a linear operator depending on $(\Psi,\Qb)$ that sends $\delta\Psi$ to a complex number \citep[cf.~e.g.,][p.~10 for differentiation in general Banach spaces]{lang1999fundamentals}.
Observe also that $h_1(t,\cdot) := \frac{\partial}{\partial t} h(t,\cdot) = \Psi_1(\cdot) - \Psi_0(\cdot)$ which is independent of $t$. 
Therefore, after stacking the complex exponentials $\{e^{i\innerprod{\kb}{\thetab}}\}_{\kb\in\Lambda_{\half}}$ into a vector $\ab(\thetab)$ in a way consistent with $\qb$, we have
\begin{equation}
	\begin{aligned}
		& J_{21} (h(t,\cdot),\Qb(t); h_1(t,\cdot)) \\
		= & - \int_{\Tbb^2} \ab(\thetab) (\Psi_t^{-1} +Q_t)^{-2} \Psi_t^{-2} (\Psi_1-\Psi_0) \d\nu_\Nb .
	\end{aligned}
\end{equation}
At this point, we are ready to convert the parametrized stationary-point equations \eqref{parametrized_eqns} to an initial value problem (IVP) of the form
\begin{equation}\label{IVP}
	\begin{cases}
		\qb'(t) = V(\qb(t),t) \\
		\qb(0) = \qb_0
	\end{cases}
\end{equation}
where the vector field
\begin{equation}
	\begin{aligned}
		V(\qb(t),t) & = -\left[ \Hcal(\Psi_t,\Qb(t)) \right]^{-1} J_{21} (h(t,\cdot),\Qb(t); h_1(t,\cdot)) \\
		& = \left[ \Hcal(\Psi_t,\Qb(t)) \right]^{-1} \times \\
		& \quad \int_{\Tbb^2} \ab(\thetab) (\Psi_t^{-1} +Q_t)^{-2} \Psi_t^{-2} (\Psi_1-\Psi_0) \d\nu_\Nb .
	\end{aligned}
\end{equation}
The following proposition summarizes the result obtained so far. Notice that here we make no distinction between the matrix $\Qb=\{q_\kb\}_{\kb\in\Lambda}$ and the ``half'' vector $\qb=\{q_\kb\}_{\kb\in\Lambda_{\half}}$.

\begin{proposition}
	Given the initial point $\qb_0$, solution to the stationary-point equation \eqref{stationary_point_eqn} with $\Psi=\Psi_0$,
	the IVP \eqref{IVP} is well-posed on the interval $t\in U$, and the terminal point $\qb(1)$ solves \eqref{stationary_point_eqn} with $\Psi=\Psi_1$.
\end{proposition}

\begin{pf}
	On the one hand, by the well-posedness results in \citet[Sec.~4]{ZFKZ2019M2}, the solution curve $\qb(t)$ is uniquely determined by the parametrized stationary-point equation \eqref{parametrized_eqns} for $t\in[0,1]$, and it also solves the IVP \eqref{IVP}.
	On the other hand, the vector field $V(\qb,t)$ on the right-hand side of the ODE in \eqref{IVP} is smooth. Therefore, the IVP is well-posed, see e.g., \citet[Appendix D]{lee2013smooth}, so its solution necessarily coincides with the desired curve $\qb(t)$ on $[0,1]$. The claim of the proposition thus follows. 
\end{pf}

In view of the proposition above, in order to solve the stationary-point equation \eqref{stationary_point_eqn} for $\Psi=\Psi_1$, one can in principle use a general-purpose ODE solver to compute the solution curve $\qb(t)$ \emph{approximately}. However, the ODE in \eqref{IVP} has a special structure due to its equivalence to the parametrized nonlinear equations \eqref{parametrized_eqns}, which enables an \emph{exact} curve-tracing algorithm called ``predictor-corrector numerical continuation''. The idea is briefly described next.

Suppose that we have reached a point on the solution curve $\qb(t)$ for some $t\in (0,1)$. From there, we aim to solve \eqref{parametrized_eqns} at $t+\delta t$ where $\delta t$ is a chosen step length. To this end, we first compute a predictor using the forward Euler method
\begin{equation}\label{Euler method}
	\begin{aligned}
		\pb(t+\delta t)=\qb(t)+V(\qb(t),t)\delta t.
	\end{aligned}
\end{equation}
Then follows a corrector step which is about solving \eqref{parametrized_eqns} at $t+\delta t$ using Newton's method initialized at the predictor $\pb(t+\delta t)$. It is hoped that Newton's method could converge to a new point $\qb(t+\delta t)$ on the solution curve. In such a way, the solution curve is continued from $t$ to $t+\delta t$, and eventually to $t=1$ by repeating such a procedure.
The main steps of the algorithm is summarized below in Algorithm~\ref{alg1:NumCont}.
\begin{algorithm}  
	\caption{Predictor-corrector numerical continuation.} 
	\label{alg1:NumCont}
	\begin{algorithmic}[1]
		\REQUIRE Initialize $t=0$, $k=0$, and $\qb^{(0)}=\textbf{0}\in \Lscr_+$. Let $\delta t=0.5$
		\WHILE{$t<1$}
		\STATE Compute a predictor $\pb^{(k+1)}$ using the Euler method \eqref{Euler method}
		\STATE Solve \eqref{parametrized_eqns} at $t+\delta t$ for $\qb^{(k+1)}$ initialized at $\pb^{(k+1)}$ using Newton's method
		\STATE Update $t:=t+\delta t$, and $k:=k+1$
		\ENDWHILE
		\ENSURE Return $\qb$ as the solution corresponding to $t=1$
	\end{algorithmic} 
\end{algorithm}

\begin{proposition}
	Algorithm~\ref{alg1:NumCont} returns the terminal point $\qb(1)$, i.e., solution to \eqref{stationary_point_eqn} with $\Psi=\Psi_1$ in a finite number of steps.
\end{proposition}

\begin{pf}
	The proof is quite standard in the literature of numerical continuation. For this reason, here we shall just outline the main idea of the proof. The reader can consult \citet[Theorem~4, p.~1095]{zhu2018wellposed} for technical details.
	
	Main idea of the proof:
	
	\begin{enumerate}
		\item It can be shown that the predictor $\pb(t+\delta t)$ always stays in the (open) feasible set \eqref{feasible_set} as long as $\delta t$ is sufficiently small.  
		\item If the predictor does not deviate too much from the solution curve $\qb(t)$, then by the famous Kantorovich theorem (see e.g., \cite{allgower2003introduction}), the standard Newton's method in each inner loop (Step 3 in Algorithm~\ref{alg1:NumCont}) is guaranteed to converge to the desired point $\qb(t+\delta t)$ on the solution curve.
		\item One can then choose a \emph{constant} step length $\delta t$, uniform in $t\in(0,1)$, so that the predictor-corrector algorithm succeeds for each step.
	\end{enumerate}
	
	Therefore, such an algorithm leads to a terminal point $\qb(1)$ within a finite number of steps, roughly equal to $1/\delta t$.
\end{pf}

\begin{remark}[Choice of $\Psi_0$]
	Typically, one takes $\Psi_0\equiv\sigma_{\zerob}$, the constant zeroth moment. In that case, the optimal spectrum is the reciprocal of a positive Laurent polynomial, known as the \emph{Maximum Entropy} (ME) solution, which can be computed analytically in the $1$-d continuous\footnote{By ``continuous'', we mean that the spectrum is supported on the whole frequency domain which is standard in spectral analysis.} setting corresponding to (nonperiodic) time series. In the current $2$-d discrete setting, the ME solution still needs to be computed in an iterative manner, but the optimal $\Qb$ is usually away from the boundary so that Newton's method is well-behaved.
\end{remark}

\section{Application to frequency estimation}\label{Sec:Fre_est}

The frequency estimation problem has numerous applications in science and engineering, e.g., signal processing in astrophysics, radars and sonars, and fault detection in rotational machines \citep{stoica2005spectral}. The problem has been extensively studied in the literature, and we shall elaborate on one solution technique via spectral analysis.

Consider a 2-d frequency estimation problem. The signal model is a superposition of complex sinusoids:
\begin{equation}\label{Vector measurements}
	y(\mathbf{t})=\sum_{j=1}^{\nu}a_{j} e^{i(\langle \boldsymbol{\theta}_{j},\mathbf{t}\rangle+\varphi_{j})}+w(\mathbf{t}), 
\end{equation}
where, 
\begin{itemize}
	\item the vector index $\tb=[t_1,t_2] \in \mathbb{Z}^2$ is such that $0\leq t_1\leq T_1$ and $0\leq t_2\leq T_2$,
	\item the integer $\nu$ denotes the number of sinusoids which is assumed to be \emph{known},
	\item $a_j \in \mathbb{R}$ is an amplitude,
	\item $\boldsymbol{\theta}_j = [{\theta_j}_1,{\theta_j}_2] \in \Tbb^2$ contains two $unknown$ normalized angular frequencies, 
	\item $\varphi_j \in \Tbb$ is an initial phase angle which is assumed to be uniformly distributed on $\left[0,2\pi\right)$ \citep[see][]{stoica2005spectral}, 
	\item and $w(\tb)$ is a complex white noise.
\end{itemize}
Taking the Fourier transform, the power spectrum 
of the random field $y(\tb)$ in \eqref{Vector measurements} is
\begin{equation}
	\Phi(\boldsymbol{\omega}) = 2\pi \sum_{j}^{\nu}a_j^2\delta(\boldsymbol{\omega}-\boldsymbol{\theta}_j)+\sigma_w^2,
\end{equation}
where $\delta\left(\cdot\right)$ is the Dirac delta and
$\sigma_w^2$ is the variance of the noise $w\left(\mathbf{t}\right)$. 
It can be seen that the spectrum is the superposition of the Dirac impulses over the frequency domain, and the location of each impulse gives an unknown frequency vector. Following this idea, we attempt to solve the estimation problem via peaks finding over the estimated spectrum.
More precisely, if we have an estimate $\hat{\Phi}$ of the underlying spectrum $\Phi$, we can take the frequency estimate 
\begin{equation}\label{peaks-finding}
	\hat{\thetab}_1=\underset{\omegab\in\Tbb^2_\Nb}{\argmax}\  \hat{\Phi}(\omegab),
\end{equation}
and $\hat{\thetab}_2$ corresponds to the second highest peak of $\hat{\Phi}$ and so on.
The general steps of the solution to the frequency estimation problem are described as follows:
\begin{enumerate}
	\item Given an index set $\Lambda$, compute the covariance lags $\{\hat\sigma_{\kb},\kb\in \Lambda\}$ of the random field \eqref{Vector measurements} via \eqref{cov_est_ave}.
	\item Solve the dual problem \eqref{J_dual} using Newton's method for an estimated power spectrum $\hat{\Phi}$ of the form \eqref{optimal solution}.
	\item Find the $\nu$ highest peaks of $\hat{\Phi}(e^{i\thetab})$ and take the corresponding arguments $\hat{\thetab}_1,\cdots,\hat{\thetab}_\nu$ as estimates of the hidden frequency vectors.
\end{enumerate}

We perform Monte-Carlo simulations with 500 repeated trials in Matlab. In generating the measurements \eqref{Vector measurements}, we have taken the integer $\nu=2$, the amplitude vector of the sinusoids $\mathbf{a}=[a_1,a_2]=[1,1]$, the noise variance $\sigma^2_w=1$, and the data size $\Tb=[T_1,T_2]=[30,30]$. In each trial, the frequency vector $\boldsymbol{\theta}_j\ (j=1,2)$ is drawn from the uniform distribution in $\Tbb^2$. In solving the dual optimization problem \eqref{J_dual}, we have taken the grid size $\Nb=\left[30,30\right]$, the index set $\Lambda$ in \eqref{set_Lambda} with the parameters $n_1=n_2=3$, and the prior $\Psi_0\equiv\sigma_{\zerob}$.

\begin{figure}
	\begin{center}
		\includegraphics[width=8.4cm]{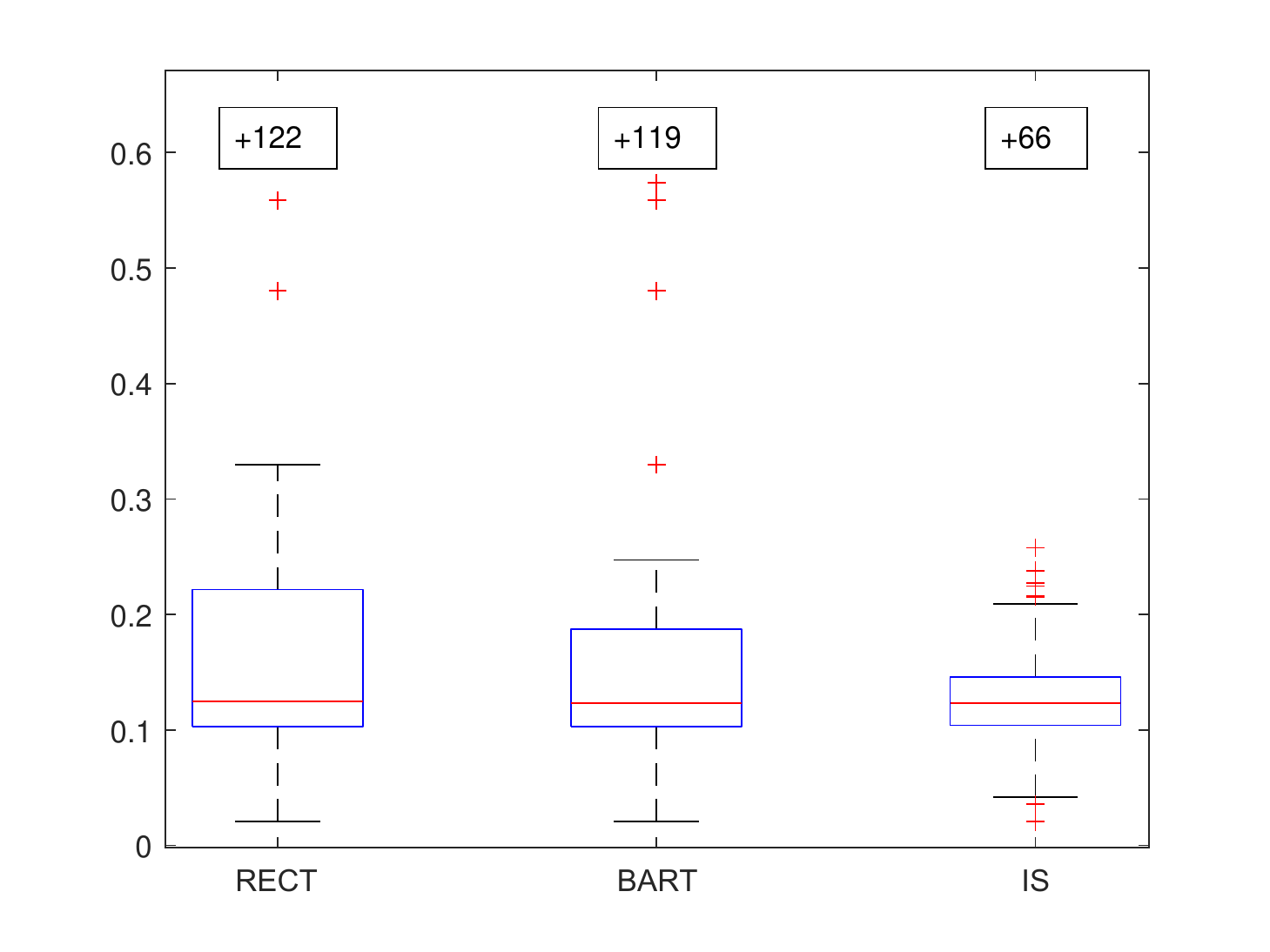}    
		\caption{Frequency estimation errors of three different spectral estimators. The central mark indicates the median, the bottom and top edges of the box indicate the 25th and 75th percentiles, respectively, and the numbers above the boxes denote the total number of outliers.}
		\label{sec2:subsec3:fg1}
	\end{center}
\end{figure}

We consider three spectral estimators in the simulations, including two windowed periodograms and our optimal IS estimator.
We utilize window functions to smooth the periodogram which is commonly done in practice \citep{engels2017advances}. More precisely, two window functions are considered: 
\begin{enumerate}
	\item the rectangular window, that is,
	\begin{equation*}
		w(\kb)=1 \ \text{ for } \ \kb\in\Lambda;
	\end{equation*}
	\item and the Bartlett window, that is,
	\begin{equation*}
		w(\kb)=w_1(k_1)w(k_2),
	\end{equation*}
	where $\kb=[k_1,k_2]\in\Lambda$, and
	\begin{equation*}
		w_j(k_j)=\frac{n_j+1-|k_j|}{n_j+1},\quad j=1,2.
	\end{equation*}
\end{enumerate}
The window widths are $\nb=[8,8]$ for the rectangular window and $\nb=[12,12]$ for the Bartlett window, and these parameters are chosen so that the windowed periodograms can exhibit good performances.
The estimates of the frequency vectors, obtained via peaks finding over each of the estimated spectra, are arranged into a long vector $\hat{\Theta}:=[{\hat{\boldsymbol{\theta}}}_1 ; {\hat{\boldsymbol{\theta}}}_2]$ where the semicolon here conforms with the Matlab syntax. The three spectral estimators are labeled as RECT, BART and IS, where RECT and BART mean that the periodogram has been computed using the rectangular window and the Bartlett window, respectively, and IS denotes our optimization approach.
We use the norm of the error defined as {$\|\hat{\Theta}-\Theta\|$} to measure the performances of the three estimators. We store the error that results in each trial, and plot the cumulative errors in each Monte-Carlo simulation using the Matlab command \verb|boxplot|. The results are depicted in Fig.~\ref{sec2:subsec3:fg1}, and one can clearly see that our IS estimator outperforms the periodogram-based methods. 
Furthermore, the number of the outliers (which the symbol $``+"$ means) of the IS approach is also clearly fewer than that of the periodograms.


In Table \ref{table_time} we show the average solution times of the three estimators. 
The IS estimator runs slower than the windowed periodograms, which is not surprising because our method requires the solution of a sophisticated optimization problem while the periodograms essentially involve only linear operations.
In plane words, we are trading the computational time for a better performance of the spectral estimator.


\begin{table}[t!]
	\begin{center}
		\caption{\centering Average solution times in 100 trails for two periodograms and the optimal IS spectrum.}	
		\label{table_time}
		\begin{tabular}{lllll}  
			\toprule   
			Methods& RECT &BART &IS \\
			\midrule   
			Avg. (s) & 0.0002064 & 0.0002175 & 0.7073   \\  
			\bottomrule  
		\end{tabular}
	\end{center}
\end{table}

\subsection{High-Resolution Property of the Spectral Estimator}

In many practical applications of the frequency estimation problem, we often face the question of resolution: when two hidden frequencies are close to each other, can they be resolved by a particular method? In the 1-d case, it is well known that the resolution limit of the periodogram is $1/T$ where $T$ is the number of measurements. In other words, when the distance between two hidden frequencies is smaller than $1/T$, the periodogram can not separate them, see e.g., \citet{BGL-THREE-00}.

In order to test the resolution property of our IS estimator, we perform three additional simulations. 
The signal model is still \eqref{Vector measurements} with $\nu=2$, and now the frequency vectors are specified in the following three cases:
\begin{enumerate}
	\item[(A)] $\thetab_1=[2.3,2.3]$ and $\thetab_2=[2.3,4.4]$,
    \item[(B)] $\thetab_1=[2.3,2.3]$ and $\thetab_2=[2.3,2.6]$,
	\item[(C)] $\thetab_1=[2.4,2.4]$ and $\thetab_2=[2.51,2.51]$. 
\end{enumerate}
The other parameters are the same as those in the previous Monte-Carlo simulations. Since the grid size in the frequency domain is $30\times30$, the minimum separation between two adjacent spectral peaks is (roughly) $2\pi/N_1=2\pi/N_2=2\pi/30 \approx 0.2094$, 
which is evidently larger than $\|\thetab_1-\thetab_2\|$ in Case C. Hence, a spectral estimator with a higher resolution is desired in the latter case in order to separate the hidden frequencies.

\begin{figure}
	\begin{center}
	   \includegraphics[width=8.4cm]{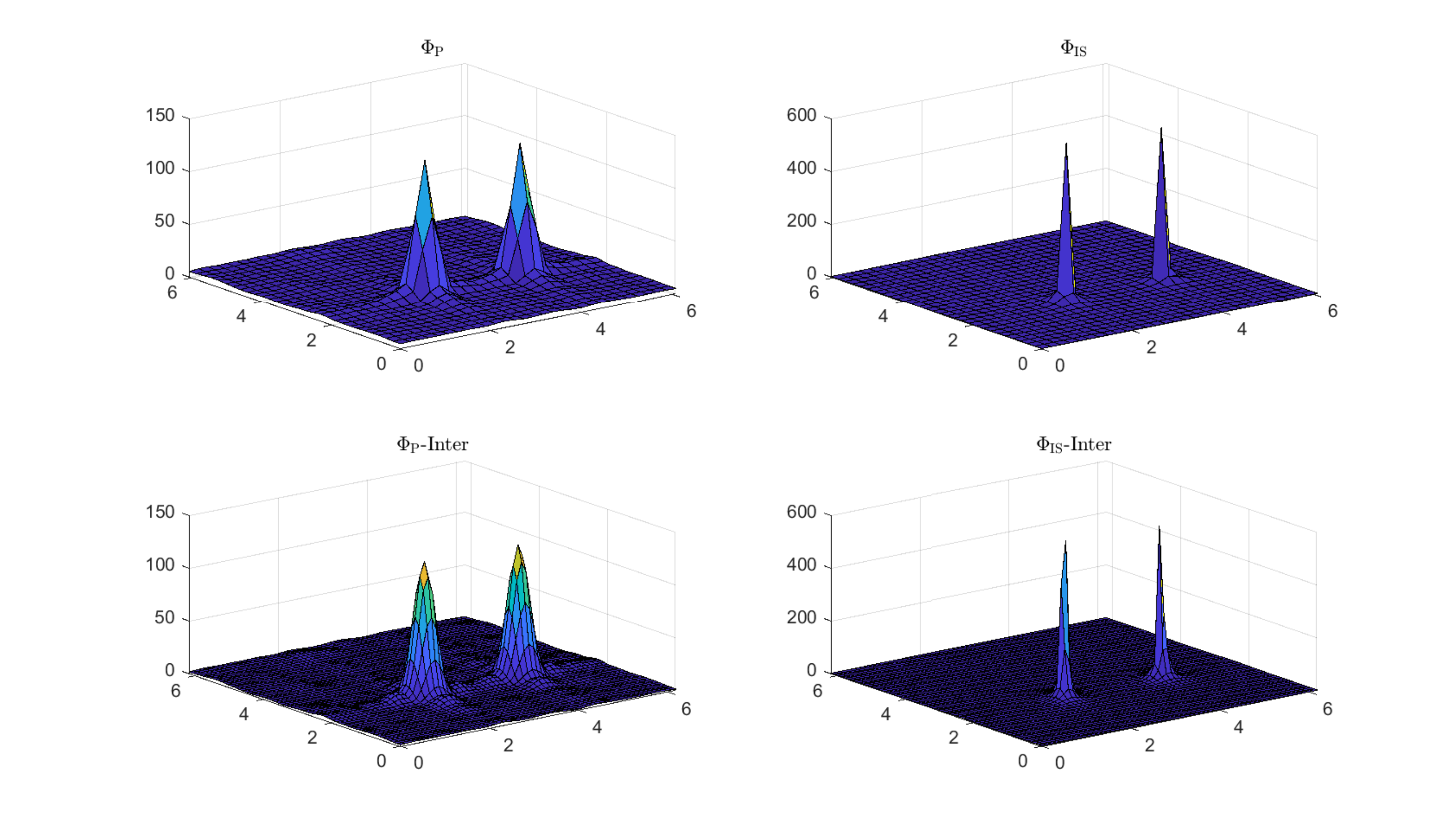}
		\caption{Case A: Estimated spectra using the windowed periodogram and the IS method, where both interpolated spectra show two proper peaks on a $60\times60$ grid. Both methods give the frequency estimates, i.e., the peak locations, $\hat{\thetab}_1=[2.3038,2.3038]$, $\hat{\thetab}_2=[2.3038,4.3982]$, the best grid points corresponding to the true frequency vector.
		}
		\label{Resolution1}
	\end{center}
\end{figure}

\begin{figure}
	\begin{center}
		\includegraphics[width=8.4cm]{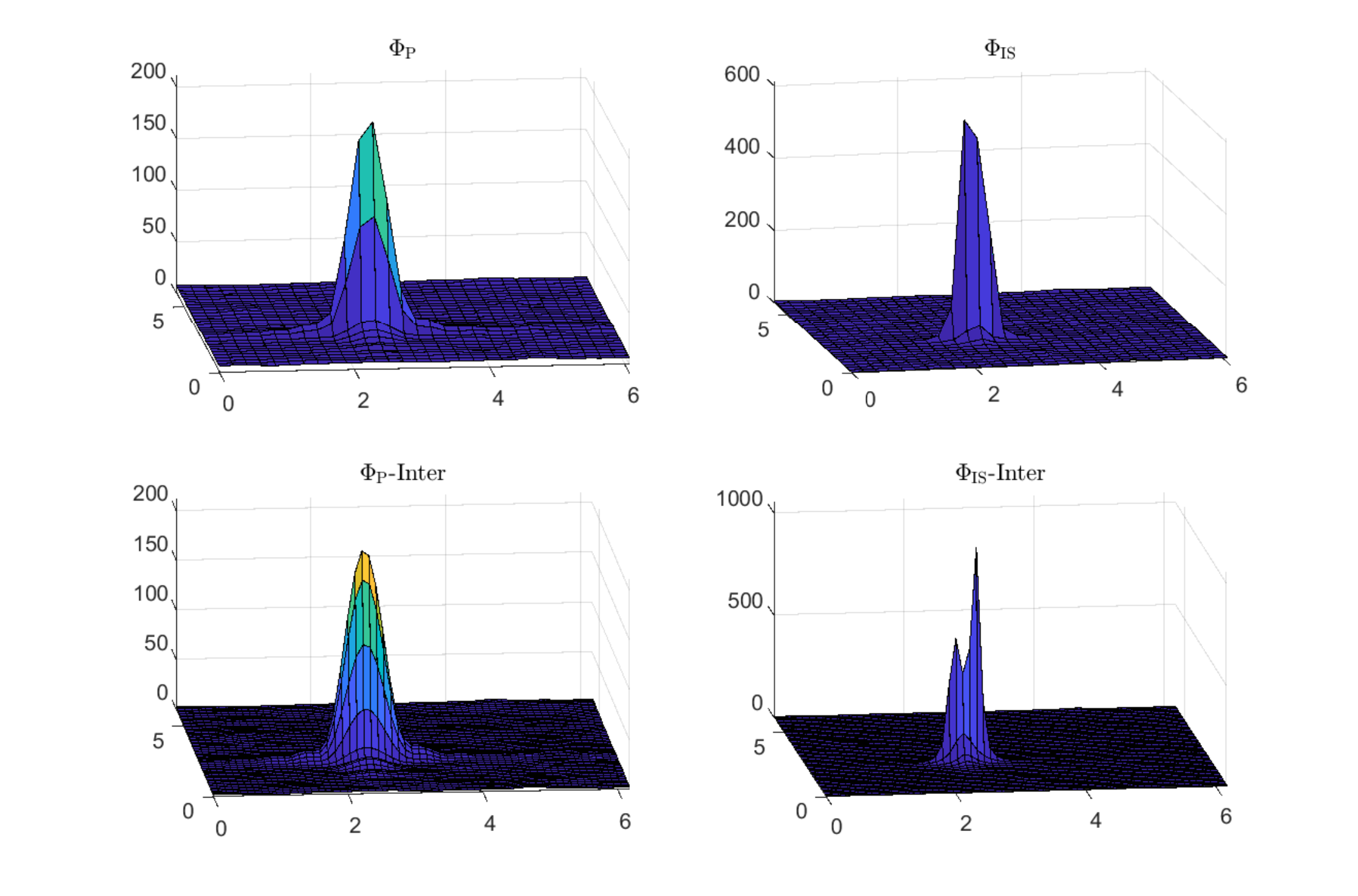}
		\caption{Case B: Estimated spectra using the windowed periodogram and the IS method, where only the IS estimator can identify two proper peaks after interpolation on a $60\times60$ grid. The IS method gives the correct the frequency estimates on the grid, i.e., the peak locations, $\hat{\thetab}_1=[2.3038,2.3038]$, $\hat{\thetab}_2=[2.3038,2.6180]$ with an error $\|\hat{\Theta}-\Theta\| = 0.0192$, while the periodogram has a large estimation error of $0.1392$.
		}
		\label{Resolution2}
	\end{center}
\end{figure}

\begin{figure}
	\begin{center}
		\includegraphics[width=8.4cm]{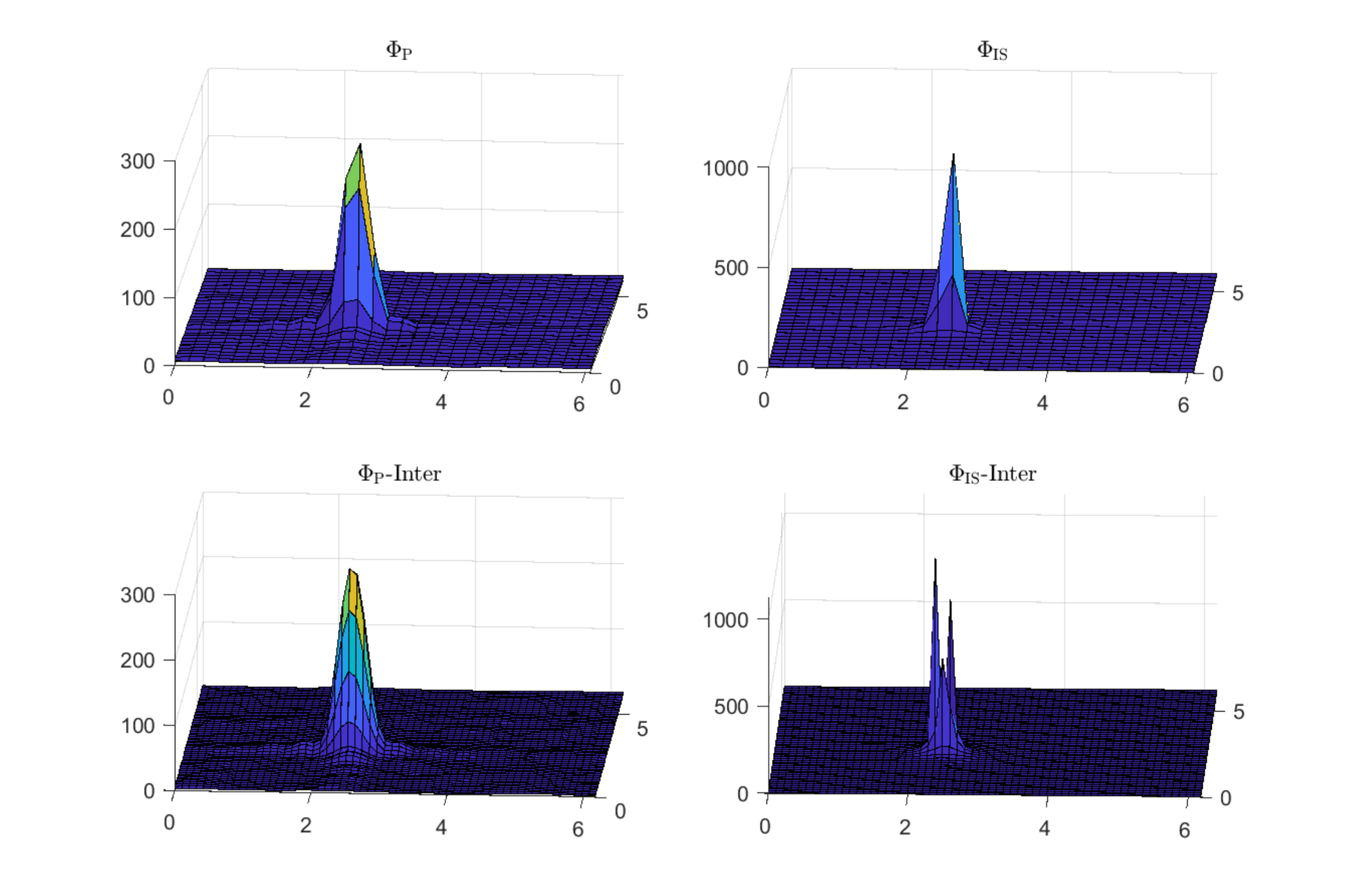}
		\caption{Case C: Estimated spectra using the windowed periodogram and the IS method, where only the IS estimator can identify two proper peaks after interpolation on a $60\times60$ grid.
        }
		\label{Resolution3}
	\end{center}
\end{figure}

The simulation results are depicted in Figs.~\ref{Resolution1}, \ref{Resolution2}, and \ref{Resolution3} corresponding to Cases A, B, and C, respectively. We consider two estimators of the spectrum $\Phi_\mathrm{P}$ and $\Phi_{\mathrm{IS}}$, where the subscript IS denotes our optimization approach and P the windowed periodogram. In view of the Monte-Carlo simulations before, here we just take the smoothed periodogram using the Bartlett window, which performs better than the other one with the rectangular window.
In each figure, ``Inter'' in the title of some subfigures means interpolating the power spectra $\Phi$ on a $60\times 60$ grid of the frequency domain,
so that peaks close to each other can possibly be separated.
	More precisely, the windowed periodogram is a trigonometric polynomial which can in principle be evaluated at any point in $\Tbb^2$.
	Similarly, in view of \eqref{optimal solution}, the optimal IS spectral density
	\begin{equation*}
		\hat{\Phi}(e^{i\thetab})=\left[(\Psi^{-1}+\hat{Q})(e^{i\thetab})\right]^{-1},
	\end{equation*}
can also be extended to any point $\thetab\in \Tbb^2$. Therefore, it is possible to evaluate both spectra on a denser grid, e.g., $60\times 60$. 

As shown in Figs.~\ref{Resolution1}, \ref{Resolution2}, and \ref{Resolution3}, the periodogram does reasonably well in Case A but fails resolving the two hidden frequency vectors in Cases B and C after interpolation, 
while the IS estimator is capable of identifying the peaks in all three cases. 
Especially in Case B, the estimation error $\|\hat{\Theta}-\Theta\|$ of the IS estimator on the $60\times60$ grid is $0.0192$. In contrast, the windowed periodogram 
locates some spurious peaks with an estimation error $0.1392$. Case C is a more challenging situation, where our IS spectrum can still exhibit two sharp peaks after interpolation.
We conclude that our IS spectral estimator not only does considerably better than traditional methods in terms of estimation errors in the Monte-Carlo simulations, but also exhibits a high-resolution property.

\subsection{Fast Inversion of the Structured Hessian}
Next we implement a numerical example on the inversion of the Hessian utilizing the its TBT structure.
More precisely, we compute one Hessian matrix $\Hcal(\Qb)$, which is positive definite of dimension $n(2n+1)+n+1$, in each trial of the Monte-Carlo simulation reported previously.
Here we take the dimensional constants $n_1=n_2=n$ in \eqref{mat_Q} for simplicity. Therefore, the matrix $C$ in \eqref{C-matrix} is an $n(2n+1)\times n(2n+1)$ TBT matrix. 
We compare two algorithms for the inversion of such $C$. One is the efficient algorithm described in \cite{wax1983efficient}, and the other is via the command \verb|^(-1)| in Matlab.

\begin{figure}[t]
	\begin{center}
		\includegraphics[width=8.4cm]{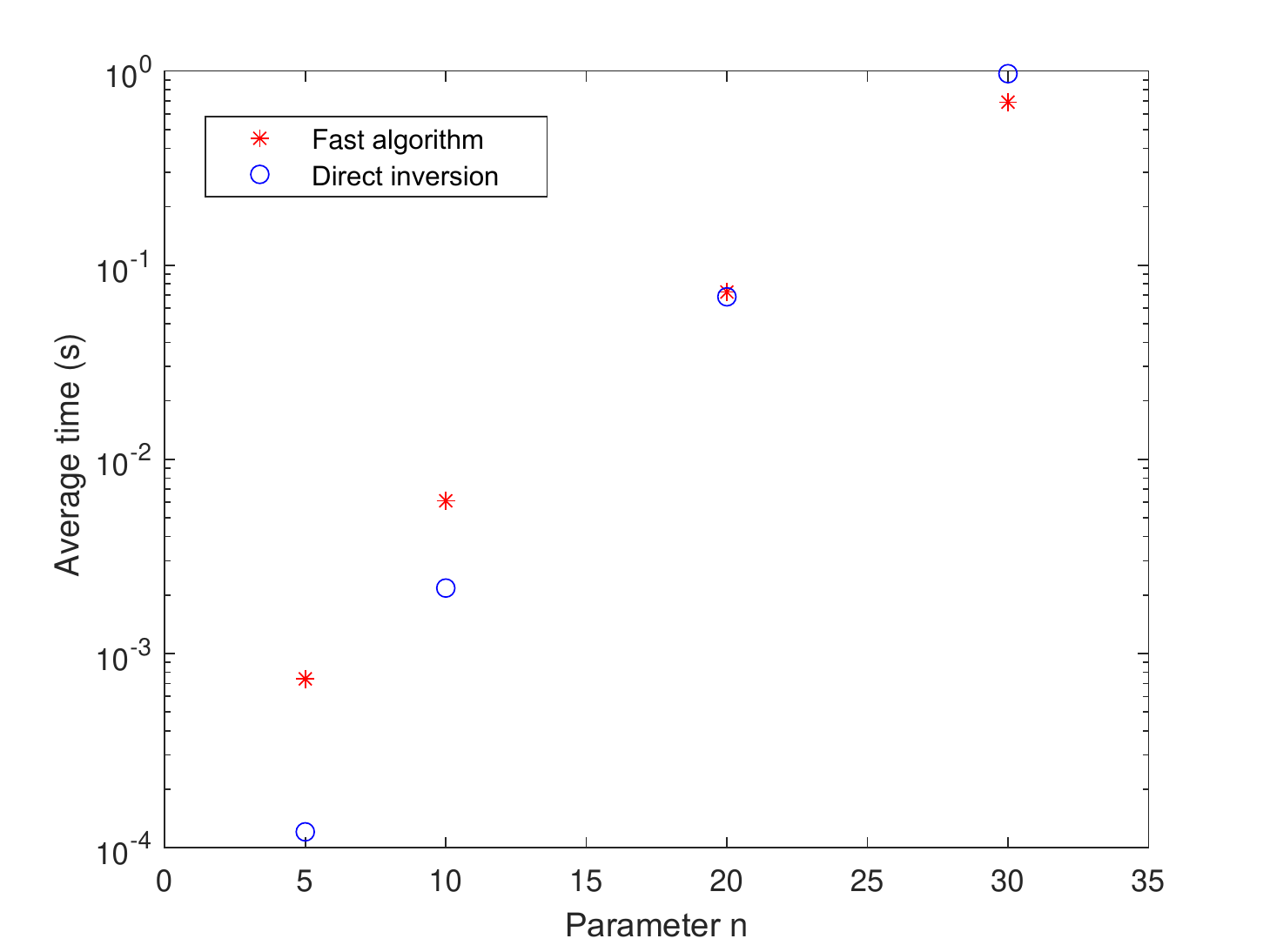}    
		\caption{Direct computation vs the fast algorithm for the inversion of the TBT matrix $C$: average solution times  with different values of the size parameter $n$.}
		\label{Average-time}
	\end{center}
\end{figure}

The results are depicted in Fig.~\ref{Average-time}. There are four groups of data points $(n,\log(t_n))$ in the figure, where $t_n$ denotes the average solution time over 100 trials for fixed $n=5,10,20$ and $30$.
One can see that when $n$ is small, a direct inversion of $C$ is less time consuming. However, when $n$ becomes large like $30$, which means that the Hessian is of a large dimension (roughly $1800$), the fast algorithm exploiting TBT structure starts to perform better. 
We must point out that such a comparison may not be fair because the fast inversion algorithm was implemented in Matlab which is known to be slower than the built-in functions. However, the fast algorithm still wins in large-scale instances, which shows its power!
In summary, when the size of the dual optimization problem is large as determined by the index set $\Lambda$, it is better to use the fast algorithm for the inversion of the Hessian in computing the Newton direction.

\section{Application to system identification}\label{sec:Sysid}

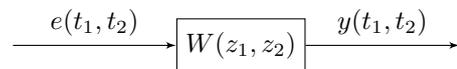
\begin{figure}[h]
	\centering
	\tikzstyle{int}=[draw, minimum size=2em]
	\tikzstyle{init} = [pin edge={to-,thin,black}]
	
	\begin{tikzpicture}[node distance=3cm,auto,>=latex']
		\node [int] (a) {$W(z_1,z_2)$};
		\node (b) [left of=a, coordinate] {};
		\node (c) [right of=a] {};
		\path[->] (b) edge node {$e(t_1,t_2)$} (a);
		\path[->] (a) edge node {$y(t_1,t_2)$} (c);
	\end{tikzpicture}
	\caption{A $2$-d linear stochastic system.}
	\label{fig:2-d_linear_system}
\end{figure}

In this section, we apply our optimization approach to solve an identification problem of 2-d linear stochastic systems, which can also be viewed as a \emph{model approximation} problem. 
Consider a 2-d linear time-invariant system as shown in Fig.~\ref{fig:2-d_linear_system}, where $W(z_1,z_2)$ is the transfer function, also called a shaping filter in signal processing \citep{SIGEST-01}. The system $W(z_1,z_2)$ is excited by a white noise process $e(t_1,t_2)$ and produces an output process $y(t_1,t_2)$. Furthermore, we assume that the transfer function has the form
\begin{equation}\label{transfer_func}
	W(\zb)=\frac{b(\zb)}{a(\zb)}=\frac{\sum_{\kb \in \Lambda_{+,2}}b_{\kb}\zb^{-\kb}}{\sum_{\kb \in \Lambda_{+,1}}a_{\kb}\zb^{-\kb}},
\end{equation}
where 
\begin{align}\label{set_Lambda_plus}
	\Lambda_{+,1} := \{(k_1,k_2)\in\Zbb^2 : 0\leq k_1,k_2 \leq n\} \\
	\Lambda_{+,2} := \{(k_1,k_2)\in\Zbb^2 : 0\leq k_1,k_2 \leq m\} \nonumber
\end{align}
with $n,m\in\Nbb_+$, $\zb=(z_1,z_2)$ is a complex vector, and $\zb^{\kb}$ stands for $z_1^{k_1}z_2^{k_2}$. 
Equivalently, the system can be described by an autoregressive moving-average (ARMA) recursion
\begin{equation}
	\sum_{\kb \in \Lambda_{+,1}} a_{\kb}\, y(\tb-\kb) = \sum_{\kb \in \Lambda_{+,2}} b_{\kb}\, e(\tb-\kb).
\end{equation}
If the white noise input is of unit variance, the power spectrum of the output process $y(\tb)$ is
\begin{equation}\label{True-spectrum}
	\Phi(e^{i\thetab})=W(e^{i\thetab})W(e^{-i\thetab})=\frac{P(e^{i\thetab})}{Q_1(e^{i\thetab})},
\end{equation}
where $P(e^{i\thetab})=|b(e^{i\thetab})|^2$ and $Q_1(e^{i\thetab})=|a(e^{i\thetab})|^2$ are Laurent trigonometric polynomials. 
In contrast, if we take the prior $\Psi=P$, the optimal form in \eqref{optimal solution} becomes
\begin{equation}\label{hat_Phi_P}
	\hat{\Phi} = (P^{-1}+\hat{Q})^{-1} = \frac{P}{1+P\hat{Q}}.
\end{equation}
Hence, $\Phi$ in \eqref{True-spectrum} and $\hat{\Phi}$ in \eqref{hat_Phi_P} specify to two model classes of different functional forms.
When the true model is of the AR type, i.e., $P\equiv1$, then clearly the two model classes coincide. 
However, when the true model is of the general ARMA type, i.e., $P$ is a nontrivial Laurent polynomial, the two model classes seem different given the observation that the denominator in \eqref{hat_Phi_P} often has a higher degree than that in \eqref{True-spectrum}.
In that case, we want to study via simulations whether it is possible to approximate the true model $\Phi$ in \eqref{True-spectrum} with some $\hat{\Phi}$ from the model class \eqref{hat_Phi_P} via matching a number of covariances. The procedure is outlined as follows:

\begin{enumerate}
	\item Given the system parameters, compute the true power spectrum $\Phi$ and the corresponding covariances $\{\sigma_{\kb}, \kb \in \Lambda\}$ where the index set $\Lambda$ is given in \eqref{set_Lambda} with $n_1=n_2=n$, the same $n$ that defines $\Lambda_{+,1}$ in \eqref{set_Lambda_plus}.
	\item Given the covariances computed in the previous step, Estimate an power spectrum $\hat{\Phi}$ of the form \eqref{optimal solution} via the solution of the dual problem \eqref{J_dual}.
\end{enumerate}

\begin{remark}
	Given the numerator polynomial $P$, a rational spectrum in the model class \eqref{True-spectrum}, more precisely the denominator $Q_1$, can be completely recovered from its covariances in the index set $\Lambda$ specified in the first step of the above procedure, by solving a suitable convex optimization problem, see \cite{RKL-16multidimensional}. The design of the procedure above for model approximation is motivated by this fact. Indeed, we keep the same prior spectral density $P$ and the same set of covariances, but solve a different convex optimization problem \eqref{optimization problem} which leads to a different model \eqref{hat_Phi_P}.
\end{remark}

For simplicity, we take $n=m=1$ in the index sets \eqref{set_Lambda_plus},
and impose a separable form $(1-\alpha_1z_{1}^{-1})(1-\alpha_2z_2^{-1})$ on the polynomials $a(\zb)$ and $b(\zb)$. 
Clearly, $\alpha_1$ and $\alpha_2$ are roots of the respective univariate polynomials, also identified as poles or zeros of the transfer function $W$. 
We take $|\alpha_j|\leq 1$ for $j=1,2$ so that the transfer function \eqref{transfer_func} is stable and minimum-phase. 
The parameters of the polynomial $a(\zb)$ can be collected into a matrix 
\begin{equation}
	\Ab=\left[\begin{matrix}
		1&-\alpha_2\\
		-\alpha_1&\alpha_1\alpha_2
	\end{matrix}\right],
\end{equation}
where $a_{k_1,k_2}=\Ab_{k_1+1,k_2+1}$ and one can define a similar matrix $\Bb$ for $b(\zb)$. 
Once the approximate spectrum $\hat\Phi$ is constructed, we compute the cumulative relative error $\|\hat{\Phi}-\Phi\|_\F/\|\Phi\|_\F$ 
where the 2-d discrete spectra can be viewed as matrices and $_{\F}$ denotes the Frobenius norm. 

In particular, we fix the $\Bb$ matrix and take three different $\Ab$ matrices as follows:
\begin{equation}\label{sys_mat_B}
	\Bb=\left[\begin{matrix}
		0.6696&-0.5357\\
		-0.4018&0.3214
	\end{matrix}\right],\ 
	\Ab_1=\left[\begin{matrix}
		1&-0.07\\
		-0.05&0.0035
	\end{matrix}\right],  
\end{equation}
\begin{equation*}
	\Ab_2=\left[\begin{matrix}
		1&-0.7\\
		-0.5&0.35
	\end{matrix}\right],\ 
	\Ab_3=\left[\begin{matrix}
		1&-0.98e^{2.1i}\\
		-0.98e^{2.1i}&0.9604e^{4.2i}
	\end{matrix}\right].
\end{equation*} 

\begin{figure}[t]
	\begin{center}
		\centerline{\includegraphics[width=1.1\columnwidth]{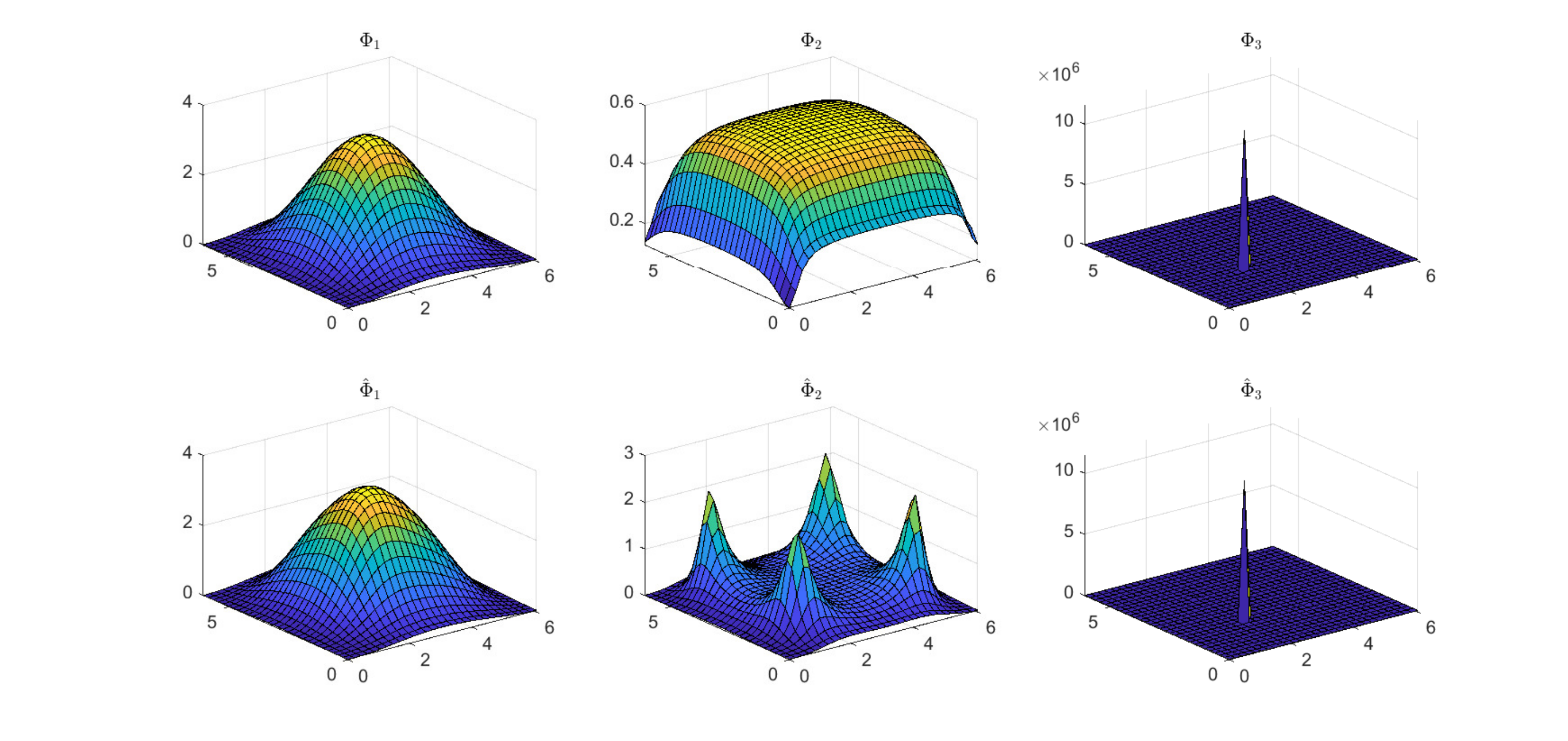}}
		\caption{The true spectra $\Phi_1$, $\Phi_2$, and $\Phi_3$ of three different systems with $\Ab_1$, $\Ab_2$, and $\Ab_3$, and their approximations $\hat{\Phi}_1$, $\hat{\Phi}_2$, and $\hat{\Phi}_3$. The relative errors of approximation for these systems are $3.25\%$, $69.88\%$, and $2.59\%$, respectively.}
		\label{Sys-iden}
	\end{center}
\end{figure}

The results of model approximation are shown in Fig.~\ref{Sys-iden}.
One can see that our method succeeds in approximating the first and third systems but fails for the second system. 
We try to explain such a success or failure next. Notice that the anti-diagonal entries of $\Ab_1$ are close to zero, meaning that the poles of the system are close to the origin, which makes the system behave like an MA one. Apparently, the rational model \eqref{hat_Phi_P} can well approximate an MA process which has a polynomial spectrum. The third system ($\Ab_3$), on the contrary, has poles very close to the unit torus, which makes the output $y(\tb)$ behave like an oscillatory signal in \eqref{Vector measurements}. This observation intuitively explains the success of model approximation in this case given the good performance of our spectral estimator in the previous section. The failure for the second system ($\Ab_2$), a generic ARMA one, is (probably) due to the difference between the model classes.

\subsection{Newton's Method vs Numerical Continuation in an Ill-Conditioned Example}\label{subsec: ill_cond_ex}

\begin{figure}[t]
	\begin{center}
		\includegraphics[width=8.4cm]{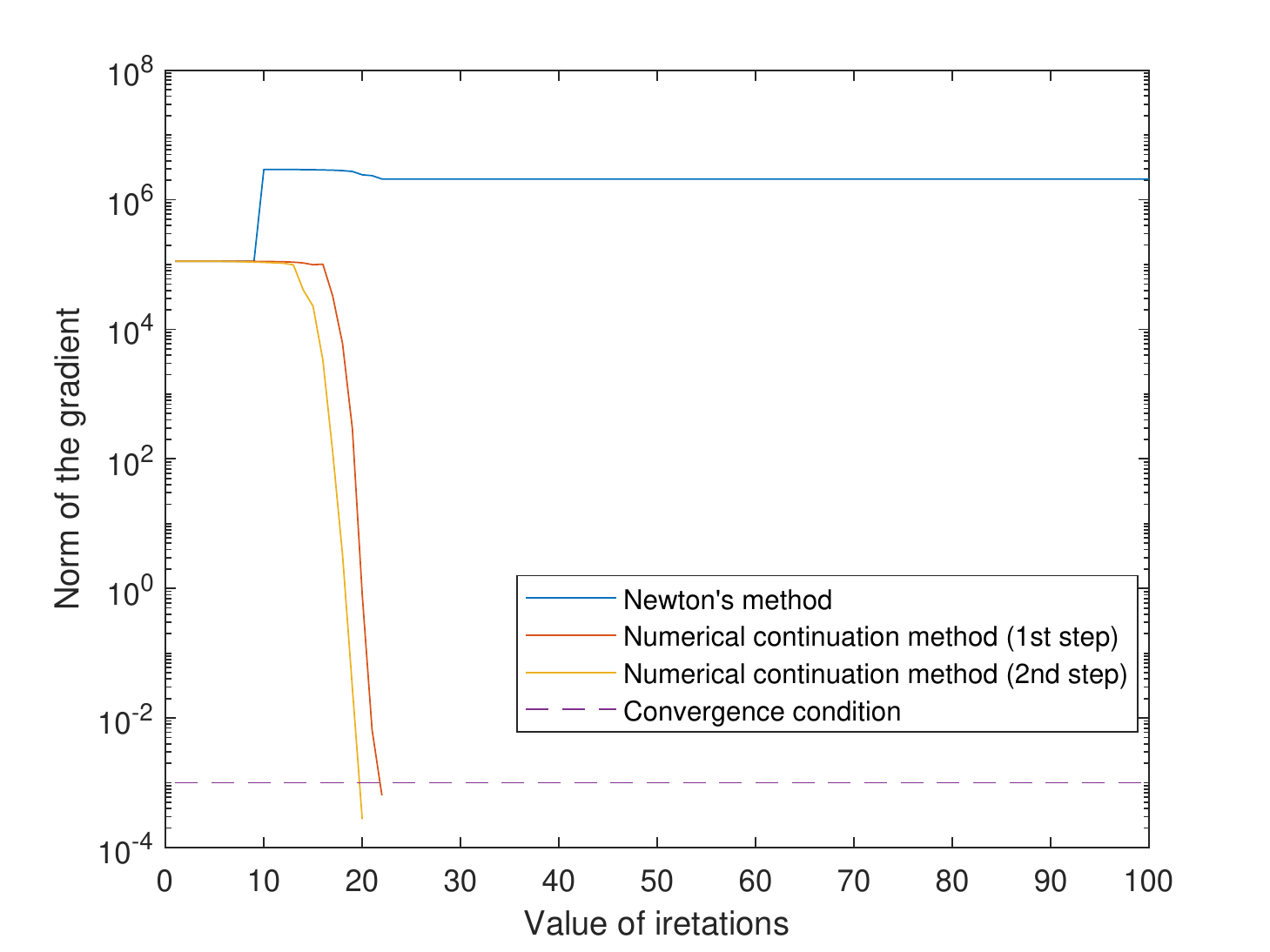}    
		\caption{Norm of the gradient changes with the number of iterations. The dashed line represents the convergence condition at the value $10^{-3}$. The blue line means that Newton's method does not converge within 100 iterations, and the other two lines show that the numerical continuation method achieves convergence for each inner loop.}
		\label{NumericalContinuationSolver}
	\end{center}
\end{figure}

In our simulations, we have found that when the poles of the true system go close to the unit torus, we face some conditioning issues. More specifically, the condition number of the Hessian becomes very large, on the scale of $10^{17}$ in several examples, which prevents the convergence of Newton's method. 
In that case, we have also observed that the dual variable in our optimization problem goes near the boundary of the feasible set.
As suggested at the beginning of Sec.~\ref{Sec:A numerical continuation solver}, 
the numerical continuation method can better handle such a situation. We illustrate this point with the following example
where we take the ``system matrix''
\begin{equation*}
	\Ab_4=\left[\begin{matrix}
		1&0.985e^{2.1i}\\
		0.985e^{2.1i}&0.9702^{4.2i}
	\end{matrix}\right],
\end{equation*} 
and keep the same $\Bb$ as in \eqref{sys_mat_B}. 
In the numerical continuation solver, we take the initial $\Psi_0 \equiv \sigma_{\zerob}$, the final $\Psi_1=P$, and the step length $\delta t=0.5$ (very large).
Then we implement our Algorithm~\ref{alg1:NumCont} in Matlab and plot how the norm of the gradient (of the dual objective function) changes with the number of iterations in Fig.~\ref{NumericalContinuationSolver}, 
where the result of Newton's method is also presented for comparison. One can see that the direct use of Newton's method performs badly (not converging), while the numerical continuation method succeeds in achieving convergence and has no issue of conditioning.
As a result, the former simply fails to approximate the underlying spectrum, while the latter succeeds with a low relative error $\|\hat{\Phi}-\Phi\|_\F/\|\Phi\|_\F=2.93\%$.

\section{Conclusions}\label{Sec:Conclusions}

We have considered a spectral estimation problem for two-dimensional second-order stationary random fields via covariance extension, i.e., 
searching a power spectrum that matches some estimated covariance lags. 
Such a problem is formulated as a convex optimization problem where the objective function is the Itakura-Saito pseudodistance between our candidate spectral density and the prior. The latter is a spectrum which embeds some \emph{a priori} information on the solution.
The corresponding dual optimization problem is solved using Newton's method, where we have exploited the fact that the Hessian of the dual objective function has a Toeplitz-block Toeplitz structure. 
In this way an efficient algorithm can be implemented for the inversion of the Hessian, which facilitate the computation of the search direction in each iteration. 
Based on the fast implementation of Newton's method, then we have developed a numerical continuation solver for our optimization problem
whose convergence is guaranteed even in some ill-conditioned instances.
We have applied our spectral estimator to problems in 2-d frequency estimation and system identification.
Simulation results for the frequency estimation problem show that the IS estimator significantly outperforms the periodogram-based estimators in terms of both the estimation errors and the frequency resolution.
In the simulations on system identification, 
we have seen that the IS estimator can well approximate a simplistic true model in the instances where the poles of the true system are either close to the origin or close to the unit torus. 
\bibliography{references}
\end{document}